\documentclass[a4paper,10pt]{article}
\usepackage{anysize}
\marginsize{2.0cm}{2.0cm}{2.0 cm}{2.0 cm}
\usepackage{setspace}

\usepackage{latexsym,amsfonts,amsmath}
\usepackage{graphicx}
\usepackage{color}
\usepackage{epstopdf}

\def\liminf{\mathop{\underline{\lim}}}
\def\limsup{\mathop{\overline{\lim}}}

\newtheorem{condition}{Condition}[section]{\bfseries}{\itshape}

{\bfseries}{\itshape}
\newtheorem{theorem}{Theorem}[section]{\bfseries}{\itshape}

{\bfseries}{\itshape}

{\bfseries}{\itshape}

{\bfseries}{\itshape}

\newtheorem{lemma}{Lemma}[section]{\bfseries}{\itshape}

{\bfseries}{\itshape}

\newtheorem{definition}{Definition}[section]{\bfseries}{\itshape}

\begin{document}

\title{Average optimality for continuous-time Markov decision processes under weak continuity conditions}
\author{Yi Zhang \thanks{Department of Mathematical Sciences, University of Liverpool, Liverpool, L69 7ZL, U.K.. E-mail: yi.zhang@liv.ac.uk and zy1985@liv.ac.uk. }}
\date{}
\maketitle \doublespace

\par\noindent{\bf Abstract:}
This article considers the average optimality for a
continuous-time Markov decision process with Borel state and
action spaces and an arbitrarily unbounded nonnegative cost rate.
The existence of a deterministic stationary optimal policy is
proved under a different and general set of conditions as compared
to the previous literature; the controlled process can be
explosive, the transition rates can be \textit{arbitrarily}
unbounded and are \textit{weakly} continuous, the multifunction
defining the admissible action spaces can be neither
compact-valued nor upper semi-continuous, and the cost rate is not
necessarily inf-compact.



\section{Introduction}

In this article we establish the existence of a deterministic
stationary average optimal policy for a possibly explosive CTMDP
(continuous-time Markov decision process) in Borel state and
action spaces under the weak continuity condition.

The average criterion for CTMDPs has been studied by many authors;
for the recent developments, see
\cite{GuoLiu:2001,GuoIEEE:2003,Guo:2006TOP,PrietoOHL:2012book} for
the case of a countable state space, and
{\cite{GuoRider:2006,GuoYe:2010AAP,Kitaev:1986,Zhu:2008} for the
case of a possibly uncountable state space. Considering a
nonnegative cost rate as in the present article, the standard
approach of proving the existence of a deterministic stationary
optimal policy for an average CTMDP is through the optimality
inequality \cite{GuoLiu:2001,GuoYe:2010AAP}. If additional but
less verifiable conditions are imposed, one can establish the
optimality equation \cite{GuoIEEE:2003,Zhu:2008}. In general, it
is known \cite{GuoLiu:2001} that the optimality equation may not
have a solution even if the optimality inequality can be solved,
see also \cite{CavazosCadena:1991}.

In the present article, for the CTMDP with Borel state and action
spaces and a nonnegative cost rate, we also follow the optimality
inequality approach, however, under the conditions different from
the present literature on CTMDPs with the average criterion. Below
we explain that our conditions are rather general, in which the
contribution of the present article also lies.

Firstly, all the aforementioned works on CTMDPs
\cite{GuoLiu:2001,GuoIEEE:2003,Guo:2006TOP,GuoRider:2006,GuoYe:2010AAP,Kitaev:1986,PrietoOHL:2012book,Zhu:2008}
assume the underlying process to be non-explosive; and most of
them achieve this by assuming the existence of a Lyapunov function
bounding the growth of the transition rates. In the present
article we remove this condition, and allow the transition rates
to be essentially arbitrarily unbounded, and the controlled
process to be possibly explosive. The development of the theory
covering such CTMDPs was once regarded quite challenging in the
survey \cite{Guo:2006TOP}; for the discounted criteria it has been
done in e.g., \cite{FeinbergBook:2012}, see also
\cite{PiunovskiyZY:2012}.

Secondly, we assume the weak continuity on the underlying signed
kernel defining the transition rates, while all the previous
literature on average CTMDPs in Borel spaces is based on the
strong continuity condition, except for \cite{GuoHuangSong:2012},
which establishes the existence of a randomized stationary optimal
policy for the constrained CTMDPs. It is relevant to point out
that recently the developments of the theory of average DTMDPs
(discrete-time Markov decision processes) and SMDPs (semi-Markov
decision processes) with weakly continuous (also called Feller)
transition probabilities have received much attention from the
research community
\cite{CostaD:2012,FeinbergLewis:2007,Feinberg:2012b,JaskiewiczN:2006JMAA,JaskiewiczN:2006OL,Jaskiewicz:2009}.
In a nutshell, as compared to the strongly continuous case, the
proofs with weakly continuous transition rates are more technical,
and the construction of the solution to the optimality inequality
would involve the notion of the generalized lower limit and the
generalized Fatou's lemma. Moreover, based on a neat
generalization of the Berge theorem \cite{FeinbergKZ:2013}, which
is partially summarized in Lemma \ref{Zy13BergeF} below, and as in
\cite{Feinberg:2012b} for the average DTMDP, we allow the
multifunction defining the admissible action spaces to be neither
compact-valued nor upper semi-continuous.

If the state space is countable, then the concepts of weak and
strong continuity coincide. However, in general, meaningful
applications of Markov control problems  to, e.g., inventory
management, have been noted, where the weak continuity condition
can be satisfied while the strong continuity condition is not, see
the examples in Section 6 of \cite{Jaskiewicz:2009}.

Since the solution to the optimality inequality is constructed
following the vanishing discount factor approach, some of the
results about discounted CTMDPs are incidentally extended in the
present paper as well.

Out of the current literature on CTMDPs, this paper is most
closely related to \cite{GuoYe:2010AAP}, which is an extension of
\cite{GuoLiu:2001}, and also derives the average optimality
inequality for a CTMDP. Nevertheless, it assumes the existence of
a Lyapunov function, and considers strongly continuous transition
rates. A more detailed comparison of our conditions with those of
\cite{GuoYe:2010AAP} is presented after Condition
\ref{Zy13GuoYeCon} below.

Finally, since we allow the transition rates to be essentially
arbitrarily unbounded and not separated from zero, the standard
technique transforming the concerned average CTMDP to an
equivalent DTMDP \cite{Puterman:1994} remains to be formally
justified and is thus not directly applicable to our setup.

The rest of this paper is organized as follows. Section
\ref{Zy13Sec2} describes the concerned CTMDP problem. The main
result is presented in Section \ref{Zy13Sec3}. The proof of the
main result is postponed to Section \ref{Zy13Sec4} with some
auxiliary statements being presented therein. We finish this
article with a conclusion in Section \ref{Zy13Sec5}. To improve
the readability, the proofs of the auxiliary results and some
definitions together with known lemmas are collected in the
appendix.

\section{Optimal control problem statement}\label{Zy13Sec2}

\par\noindent\textbf{Notations and conventions.} In what follows, $I$ stands for the indicator function, $\delta_{x}(\cdot)$
is the Dirac measure concentrated at $x,$ and ${\cal{B}}(X)$ is
the Borel $\sigma$-algebra of the topological space $X.$ Below,
unless stated otherwise, the term of measurability is always
understood in the Borel sense, and a function can take values in
$[-\infty,\infty].$ The convention of $\infty-\infty:=\infty$ is
in use.

The primitives of a CTMDP are the following elements $\{S, A,
(A(x)\subseteq A, x\in S), q(\cdot|x,a)\},$ where $S$ is a
nonempty Borel state space, i.e., a measurable subset of some
complete separable metric space, $A$ is a nonempty Borel action
space, and the multifunction $A(\cdot):$ $x\mapsto A(x)\subseteq
A$ specifies the admissible action spaces, for which we assume
that $A(x)\in {\cal B}(A)$ for each $x\in S$, and its graph
$K:=\{(x,a):x\in S,a\in A(x)\}$ belongs to ${\cal B}(S\times A)$
and contains the graph of at least one measurable mapping from $S$
to $A$. This assumption guarantees the existence of deterministic
stationary policies defined below. The transition rates are given
by $q(\cdot|x,a),$ a signed kernel on ${\cal{B}}(S)$ given
$(x,a)\in K$ such that $q(\Gamma_S\setminus\{x\}|x,a)\ge 0$ for
all $\Gamma_S\in{\cal{B}}(S).$ Throughout this article we assume
that $q(\cdot|x,a)$ is conservative and stable, i.e., $q(S|x,a)=0$
and $\bar{q}_x=\sup_{a\in A(x)}q_x(a)<\infty,$ where $q_x(a):=
-q(\{x\}|x,a).$

Following the Kitaev construction of a CTMDP \cite{Kitaev:1986},
we take the sample space $\Omega:= S\times((0,\infty]\times
S_\infty)^\infty$, where $S_{\infty}:=S\bigcup\{x_\infty\}$ with
the isolated point $x_{\infty}\notin S$.  We equip $\Omega$ with
its Borel $\sigma$-algebra $\cal F$. For each $n\geq 0$, and any
element $\omega:=(x_0,\theta_1,x_1,\theta_2,\dots)\in \Omega$, let
$ t_n(\omega):=t_{n-1}(\omega)+\theta_n$ with $t_0(\omega):=0,$
and $ t_\infty(\omega):=\lim_{n\rightarrow\infty}t_n(\omega).$
Obviously, $t_n(\omega)$ are measurable mappings on the sample
space $\Omega$. In what follows, we will omit the argument
$\omega\in \Omega$ from the presentation for simplicity, and
understand $t_n,$ $x_n$, $\theta_{n+1}$, and $t_\infty$ as the
$n$-th jump moment, jumpped-in state, holding time of $x_n$, and
the explosion moment. The pairs $\{t_n,x_n\}$ form a marked point
process with the internal history $\{{\cal F}_t\}_{t\ge 0}$ (see
Chapter 4 of \cite{Kitaev:1995}), which defines the stochastic
process on $(\Omega,{\cal F})$ of interest $\{\xi_t,t\ge 0\}$ by
\begin{eqnarray}\label{GZCTMDPdefxit}
\xi_t=\sum_{n\ge 0}I\{t_n\le t<t_{n+1}\}x_n+I\{t_\infty\le
t\}x_\infty,
\end{eqnarray}
where $x_\infty$ is the cemetery point so that
$A(x_\infty):=\{a_\infty\}$ and $q_{x_\infty}(a_\infty):=0$ with
$a_\infty\notin A$ being some isolated point. Below we denote
$A_\infty:=A\bigcup\{a_\infty\}.$ As in \cite{Gihman:1975} we
formally put $\xi_\infty:=x_\infty.$

\begin{definition}
A (randomized history-dependent) policy $\pi$ for the CTMDP is
given by a sequence $(\pi_n)$ such that, for each $n=0,1,\dots,$
$\pi_n(da|x_0,\theta_1,\dots,x_n,s)$ is a stochastic kernel on $A$
concentrated on $A(x_n)$, and for each
$\omega=(x_0,\theta_1,x_1,\theta_2,\dots)\in \Omega$, $t> 0,$
\begin{eqnarray*}
\pi(da|\omega,t):=I\{t\ge t_\infty\}\delta_{a_\infty}(da)+
\sum_{n=0}^\infty I\{t_n< t\le
t_{n+1}\}\pi_n(da|x_0,\theta_1,\dots,{x_n},t-t_n).
\end{eqnarray*}
\end{definition}
In other words, a policy $\pi$ is a predictable (with respect to
$\{{\cal F}_{t}\}_{t\ge 0}$) stochastic kernel  from
$\Omega\times(0,\infty)$ to $A_\infty,$ see Theorem 4.19 in
\cite{Kitaev:1995}. The class of all policies for the CTMDP is
denoted by $\Pi.$ A policy is called Markov if it is in the form
$\pi(da|\omega, t)=\pi(da|\xi_{t-}(\omega),t),$ where, with
conventional abuse of notations, $\pi$ on the right hand is a
stochastic kernel. Denote by $\Pi_M\subset\Pi$ the set of Markov
policies.

Under a policy $\pi:=(\pi_n)\in \Pi$, we define the following
random measure on $S\times (0,\infty)$
\begin{eqnarray*}
\nu^\pi(dt, dy)&:=& \int_A q(dy\setminus\{\xi_{t_-}(\omega)\}|\xi_{t-}(\omega),a)\pi(da|\omega,t)dt\nonumber\\
&=&\sum_{n\ge 0}\int_A
q(dy\setminus\{x_n\}|x_n,a)\pi_n(da|x_0,\theta_1,\dots,x_n,t-t_n)I\{t_n<
t\le t_{n+1}\}dt
\end{eqnarray*}
with $ q(dy|x_\infty,a_\infty):=0. $ Suppose that an  initial
distribution $\gamma$ on $S$ is given. Then by Theorem 4.27 in
\cite{Kitaev:1995}, there exists a unique probability measure
${P}^\pi_\gamma$ such that
\begin{eqnarray*}
{P}_{\gamma}^\pi(\xi_0\in dx)=\gamma(dx),
\end{eqnarray*}
and with respect to $P_\gamma^\pi,$ $\nu^\pi$ is the dual
predictable projection of the random measure of the marked point
process $\{t_n,x_n\}.$ The process $\{\xi_t\}$ defined by
(\ref{GZCTMDPdefxit}) under the probability measure
${}{P}_\gamma^\pi$ is called a CTMDP. Below, when $\gamma(\cdot)$
is a Dirac measure concentrated at $x\in S,$ we use the denotation
${}{P}_x^\pi.$ Expectations with respect to ${}{P}_\gamma^\pi$ and
${}{P}_x^\pi$ are denoted as ${}{E}_{\gamma}^\pi$ and
${}{E}_{x}^\pi,$ respectively. In fact, in what follows, we often
write $P^\pi$ instead of $P_\gamma^\pi$ when there is no
confusion. Under the probability measure $P_\gamma^\pi,$ the
system dynamics of a CTMDP can be described as follows. The
initial state $x_0$ has the distribution given by $\gamma.$ Given
the current state $x_{n}$, the sojourn time $\theta_{n+1}$ has the
tail function given by $ P^\pi(\theta_{n+1}\ge t |
x_0,\theta_1,\dots,x_n)=e^{-\int_0^t
\int_{A}q_{x_n}(a)\pi_{n}(da|x_0,\theta_1,\dots,x_n,s)ds}, $ and
upon a jump, the distribution of the next state $x_{n+1}$ is given
by $ P^\pi(x_{n+1}\in
\Gamma|x_0,\theta_1,\dots,x_n,\theta_{n+1})=\frac{\int_A
q(\Gamma\setminus\{x_n\}|x_n,a)\pi_{n}(da|x_0,\theta_1,\dots,x_n,\theta_{n+1})}{\int_A
q_{x_n}(a)\pi_{n}(da| x_0,\theta_1,\dots,x_n,\theta_{n+1})} $ for
each $\Gamma\in {\cal B}(S),$ where and below we quite formally
put $\int_A
q(\Gamma\setminus\{x_n\}|a)\pi_{n}(da|x_0,\theta_1,\dots,x_n,\infty):=0$
for each $\Gamma\in {\cal B}(S)$ and use the convention of
$\frac{0}{0}:=0,$ so that $
P^\pi(x_{n+1}=x_\infty|x_0,\theta_1,\dots,x_n,\theta_{n+1})=1-P^\pi(x_{n+1}\in
S|x_0,\theta_1,\dots,x_n,\theta_{n+1}). $ According to
\cite{FeinbergS:2013}, under each Markov policy $\pi$, the process
$\xi_t$ is a Markov jump process in the sense of
\cite{Gihman:1975} with respect to $(\Omega, {\cal F}, \{{\cal
F}_t\}_{t\ge 0}, P_x^\pi)$ for each $x\in S.$

We are also interested in policies in more specific forms.
\begin{definition}
With slight but conventional abuse of denotations, a policy
$\pi=(\pi_n)_{n=0,1,\dots}\in \Pi$  is called (randomized)
stationary if each of the stochastic kernels $\pi_n$ reads
$\pi_n(da|x_0,\theta_1,\dots,x_n,t-t_n)=\pi(da|x_n)$. A stationary
policy is further called deterministic if
$\pi_n(da|x_0,\theta_1,\dots,x_n,t-t_n)=\delta_{\varphi(x_n)}(da)$
for some measurable mapping $\varphi$ from $S$ to $A$ such that
$\varphi(x)\in A(x)$ for each $x\in S$; the existence of such a
mapping is guaranteed by the assumption imposed on the
multifunction $A(\cdot),$ which also implies the set $\Pi$ being
nonempty.
\end{definition}

Let $c(x,a),$ a measurable function on $K$ that takes values in
$[0,\infty)$, represent the cost rate at the present state $x\in
S$ and action $a\in A(x)$. Quite formally, for any measurable
function $f$ on $K$, we put $f(x_\infty,a_\infty)=0$. This
agreement, together with (\ref{GZCTMDPdefxit}) and that
$q(\{x_\infty\}|x_\infty,a_\infty)=0=q_{x_\infty}(a_\infty)$,
allows one to define formally the long-run average cost by
\begin{eqnarray*}
W(x,\pi)&:=& \limsup_{T\rightarrow\infty}\frac{1}{T}{E}_x^\pi\left[\int_0^{T} \int_A c(\xi_t,a)\pi(da|\omega,t)dt \right] \nonumber \\
&=&
\limsup_{T\rightarrow\infty}\frac{1}{T}{E}_x^\pi\left[\int_0^{\min\{T,~t_\infty\}}
\int_A c(\xi_t,a)\pi(da|\omega,t)dt \right]. \label{P1a}
\end{eqnarray*}
We are interested in the following optimal control problem
\begin{eqnarray}\label{Zy13probAve}
W(x,\pi) \rightarrow \min_{\pi\in \Pi},~x\in S,
\end{eqnarray}
for which a policy $\pi^\ast$ is called optimal if
$W(x,\pi^\ast)=\inf_{\pi\in \Pi}W(x,\pi)$ for each $x\in S.$


The objective of the present article is to show the existence of a
deterministic stationary optimal policy under the weak continuity
conditions on the transition rates, which can be essentially
arbitrarily unbounded.

\section{Main result}\label{Zy13Sec3}
\begin{condition}\label{Zy13G}
 $\inf_{x\in S}\inf_{\pi\in
\Pi}W(x,\pi)<\infty.$
\end{condition}

For each real constant $\alpha>0$, we define the expected total
discounted cost under each policy $\pi\in \Pi$ by
\begin{eqnarray*}
W_\alpha(x,\pi):=E_x^\pi\left[\int_0^\infty e^{-\alpha t}\int_A
c(\xi_t,a)\pi(da|\omega,t)dt\right]=E_x^\pi\left[\int_0^{\min\{t_\infty,\infty\}}
e^{-\alpha t}\int_A c(\xi_t,a)\pi(da|\omega,t)dt\right],
\end{eqnarray*}
and the value function for the corresponding discounted problem by
$W_\alpha(x):=\inf_{\pi\in \Pi}W_\alpha(x,\pi).$  Let
\begin{eqnarray*}
m_\alpha:=\inf_{x\in S}W_\alpha(x)~\mbox{and }
h_\alpha(x):=W_\alpha(x)-m_\alpha\ge 0,
\end{eqnarray*}
where the regulation of $\infty-\infty:=\infty$ is in use. The
function $h_\alpha$ on $S$ is sometimes called the relative
difference or normalized value function for the discounted
problem, on which we impose the following condition, where $\rho$
denotes the predetermined metric on $S$ consistent with its
topology.
\begin{condition}\label{Zy13ConRelative}
$\liminf_{0<\alpha\downarrow 0,y\rightarrow
x}h_\alpha(y):=\sup_{\delta>0,~\Delta>0}\left\{\inf_{0<\alpha\le
\delta,~\rho(x,y)< \Delta}h_\alpha(y)\right\}<\infty$ for each
$x\in S.$
\end{condition}
What was defined in the above condition is the generalized lower
limit of the function $h_\alpha(y)$ as $0<\alpha\downarrow 0$ and
$y\rightarrow x$. Condition \ref{Zy13ConRelative} is equivalent to
that for each $x\in S,$ there exist sequences
$0<\alpha_n\downarrow 0$ and $y_n\rightarrow x$ such that
$\{h_{\alpha_n}(y_n)\}$ is bounded. Condition
\ref{Zy13ConRelative} and its synonyms are widely assumed in the
current literature on average CTMDPs. We provide more insights on
Condition \ref{Zy13ConRelative} after we introduce Condition
\ref{Zy13WeakCon} below.

Finally, we assume the following weak continuity condition. To
this end, we recall that a function $c$ on the space $K$ is called
$\mathbb{K}$-inf-compact if it is lower semi-continuous on $K$,
and satisfies the following; for each $S\ni x_n\rightarrow x\in S$
as $n\rightarrow\infty$, each sequence $a_n\in A(x_n)$ such that
$c(x_n,a_n)$ is bounded from the above, admits a limit point $a\in
A(x)$ \cite{FeinbergKZ:2013}. The function $c$ is called
inf-compact on $K$ if the set $\{(x,a)\in K: c(x,a)\le \lambda \}$
is compact in $K$ for each $\lambda\in (-\infty,\infty).$ By the
way, the inf-compactness on $K$ is defined in a weaker sense in
\cite{Hernandez-Lerma:1996}. It is known that the inf-compactness
of a function implies its $\mathbb{K}$-inf-compactness
\cite{FeinbergKZ:2013}.

\begin{condition}\label{Zy13WeakCon}
\par\noindent(a) For each bounded continuous function $f$ on $S$,
$\int_{S} f(y)q(dy|x,a)$ is continuous in $(x,a)\in K.$
\par\noindent(b) The cost rate $c$ is $\mathbb{K}$-inf-compact.
\par\noindent(c) There exists a continuous function
$w$ on $S$ taking values in $(0,\infty)$ such that
$\overline{q}_x\le w(x)$ for each $x\in S.$
\end{condition}
The rather weak part (c) of the previous condition is for
technical convenience, and essentially allows the transition rates
to be arbitrarily unbounded, since so can be the function $w$.
Part (a) of Condition \ref{Zy13WeakCon} reads that the transition
rates are weakly continuous. Condition \ref{Zy13WeakCon} does not
require the multifunction $A(x)$ to be either compact-valued or
upper semi-continuous.

Some comments on Condition \ref{Zy13ConRelative} are in position
now. Suppose Conditions \ref{Zy13G} and \ref{Zy13WeakCon} are
satisfied, so that for any $\alpha>0,$ there exists a
deterministic stationary optimal policy $\varphi_\alpha^\ast$ for
the discounted problem, i.e.,
$W_\alpha(x)=W_\alpha(x,\varphi^\ast_\alpha)$ for each $x\in S;$
and for all sufficiently small $\alpha>0,$ $m_\alpha<\infty$ (as
explained in the proof of Theorem \ref{Zy13MainThm} below). Assume
that there exists some $z\in S$ such that for all sufficiently
small $\alpha>0,$ $m_\alpha=W_\alpha(z)$. (In fact, if
$S=\{0,1,2,\dots\}$ or $S=[a,b)$ with $a\in \mathbb{R}$ and $b\in
\mathbb{R}\bigcup \{+\infty\},$ then this assumption is satisfied
when $A(x)$ is decreasing in $x\in S,$ and for all sufficiently
small $\alpha>0,$ $\frac{c(x,a)}{\alpha+w(x)}$ and
$\frac{w(x)}{\alpha+w(x)}\int_S u(y)
\left(\frac{q(dy|x,a)}{w(x)}+I\{x\in dy\} \right)$ are increasing
in $x\in S$ for each fixed $a\in A(x)$ and increasing nonnegative
function $u$ on $S.$ This follows from the fact that
$W_\alpha(x)=\lim_{n\uparrow \infty}v_n(x)$ with $v_n$ being
defined in the proof of Lemma \ref{Zy13discountedLem} below.)
Consider the stopping time $\tau_z=\inf\{t\ge 0:\xi_t=z\}$ (with
respect to $\{{\cal F}_t\}_{t\ge 0}$). As usual, the infimum taken
over the empty set is put as $+\infty.$ It is known
\cite{Kitaev:1995} that $\min\{\tau_z,t_\infty\}$ is also a
stopping time. Then for all sufficiently small $\alpha>0,$
\begin{eqnarray*}
h_\alpha(x) &=&
E_x^{\varphi^\ast_\alpha}\left[\int_0^{\min\{\tau_z,t_\infty\}}
e^{-\alpha t}c(\xi_t,\varphi^\ast(\xi_t))dt \right] +
E_x^{\varphi^\ast_\alpha}
\left.\left[E_x^{\varphi^\ast_\alpha}\left[\int_{\min\{\tau_z,t_\infty\}}^{t_\infty}
e^{-\alpha t}c(\xi_t,\varphi^\ast(\xi_t))dt \right|{\cal
F}_{\min\{\tau_z,t_\infty\}} \right]\right]\\
&&- W_\alpha(z).
\end{eqnarray*}
Furthermore, by Theorem 4 on p.197 of \cite{Gihman:1975} the
process $\xi_t$ is a strong Markov one with respect to $\{{\cal
F}_t\}_{t\ge 0}.$ So by applying the strong Markov property to the
second summand on the right hand side of the previous equality, we
see
\begin{eqnarray*}
h_\alpha(x)&\le&
E_x^{\varphi^\ast_\alpha}\left[\int_0^{\min\{\tau_z,t_\infty\}}
e^{-\alpha t}c(\xi_t,\varphi^\ast(\xi_t))dt \right] +
E_x^{\varphi^\ast_\alpha}\left[e^{-\alpha \min\{\tau_z,t_\infty\}}
W_\alpha(z)
\right]- W_\alpha(z)\\
&\le&
E_x^{\varphi^\ast_\alpha}\left[\int_0^{\min\{\tau_z,t_\infty\}}
e^{-\alpha t}c(\xi_t,\varphi_\alpha^\ast(\xi_t))dt \right]
\le\sup_{\pi} E_x^{\pi}\left[\int_0^{\min\{\tau_z,t_\infty\}}
\int_A c(\xi_t,a)\pi(da|\omega,t))dt \right]
\end{eqnarray*}
where the first inequality further follows from the fact that
$\xi_t$ is right-continuous and $W_\alpha(x_\infty)=0.$  It can be
shown \cite{Bertsekas:1978,GuoZhang:2013} that if there is some
constant $\epsilon>0$ such that $q_x(a)>\epsilon$ for all $x\ne z$
and $a\in A(x),$ then
\begin{eqnarray}\label{ZyVerifyCon2}
\sup_{\pi} E_x^{\pi}\left[\int_0^{\min\{\tau_z,t_\infty\}} \int_A
c(\xi_t,a)\pi(da|\omega,t))dt \right]<\infty
\end{eqnarray}
for each $x\in S$ if there exists a real-valued upper
semi-analytic function $v$ on $S$ such that
\begin{eqnarray*}
0\ge c(x,a)+\int_{S\setminus \{z\}}q(dy|x,a)v(y)
\end{eqnarray*}
for each $x\ne z$ and $a\in A(x).$ This thus provides a sufficient
condition imposed on the primitives of the CTMDP model for
verifying Condition \ref{Zy13ConRelative}, which does not refer to
the existence of a Lyapunov function as in Condition
\ref{Zy13GuoYeCon} below (cf. \cite{GuoYe:2010AAP}). By the way,
if the process $\xi_t$ is non-explosive as prevailingly assumed in
the current literature, then (\ref{ZyVerifyCon2}) is satisfied
when, for example, the process $\xi_t$ exhibits some version of
the ergodic property.

Similar versions of Conditions \ref{Zy13G}, \ref{Zy13ConRelative}
and parts (a,b) of the previous condition are assumed in
\cite{Feinberg:2012b} but for discrete-time problems, see
Assumptions G, W* and \underline{B} therein.

\begin{theorem}\label{Zy13MainThm}
Suppose Conditions \ref{Zy13G}, \ref{Zy13ConRelative} and
\ref{Zy13WeakCon} are satisfied. Then there exist a constant $g$,
a nonnegative real-valued lower semi-continuous function $h$
on $S$ and a deterministic stationary policy $\varphi^\ast$
such that
\par\noindent(a) the following optimality inequality is satisfied for each $x\in S$
\begin{eqnarray} \label{Zy13OptIne}
g+w(x)h(x)&\ge& \inf_{a\in A(x)}\left\{c(x,a)+w(x)\int_S
h(y)\left(\frac{q(dy|x,a)}{w(x)}+I\{x\in
dy\}\right)\right\}\nonumber\\
&=&
c(x,\varphi^\ast(x))+w(x)\int_S h(y)\left(\frac{q(dy|x,\varphi^\ast(x))}{w(x)}+I\{x\in
dy\}\right);
\end{eqnarray}
\par\noindent(b) the deterministic stationary policy $\varphi^\ast$ is optimal for the average CTMDP problem
(\ref{Zy13probAve}); and
\par\noindent(c)
$g=\inf_{\pi\in \Pi}W(x,\pi)<\infty$ for each $x\in S.$
\end{theorem}
The proof of this theorem is postponed to the next section, by
inspecting which one can see that any deterministic stationary
policy that satisfies (\ref{Zy13OptIne}) is optimal. Furthermore,
it follows from Lemma \ref{Zy13GYLem} below that $\inf_\pi
W(x,\pi)$ is given by the smallest constant $g$ satisfying the
inequality (\ref{Zy13OptIne}).

The statement of Theorem \ref{Zy13MainThm} is obtained in
\cite{GuoYe:2010AAP} under the following Condition
\ref{Zy13GuoYeCon}, see Assumptions A, B and C therein.

\begin{condition}\label{Zy13GuoYeCon}
\par\noindent(a) There exists a measurable function $w\ge 1$ on $S$ and
constants $c_0\in (-\infty,\infty)$, $b_0\ge 0$ and $M_0> 0$ such
that
\par\noindent (i) $\int_S w(y)q(dy|x,a)\le c_0w(x)+b_0$ for
each $(x,a)\in K;$ and
\par\noindent (ii) $\overline{q}_x\le
M_0w(x)$ for all $x\in S.$
\par\noindent(b) For some sequence $\alpha_n\downarrow 0$ as
$n\uparrow\infty$ and some fixed $x_0\in S,$ there exist a real
constant $L^\ast$ and a finitely valued nonnegative measurable
function $U$ on $S$ such that
\par\noindent  (i) $\sup_{n=1,2,\dots}\left\{\alpha_n
W_{\alpha_n}(x)\right\}<\infty$ for each $x\in S;$ and
\par\noindent (ii) $L^\ast\le W_{\alpha_n}(x)-W_{\alpha_n}(x_0)\le
U(x)<\infty$ for each $x\in S.$
\par\noindent(c) The following compactness-continuity condition is satisfied.
\par\noindent (i) The set $A(x)$ is compact for each $x\in S$;
\par\noindent (ii) the cost rate $c(x,a)$ is lower semi-continuous in $a\in A(x)$
for each $x\in S;$ and
\par\noindent (iii) for each bounded
measurable function $f$ on $S$, $\int_S f(y)q(dy|x,a)$ is
continuous in $a\in A(x)$ for each $x\in S.$
\end{condition}
The function $w$ in part (a) of the above condition is called a
Lyapunov function or a bounding function, whose existence
guarantees the process $\xi_t$ to be non-explosive, i.e.,
$P_x^\pi(t_\infty=\infty)=1$ for each $x\in S$
\cite{GuoYe:2010AAP}, which is also prevailingly assumed in the
previous literature on CTMDPs with possibly unbounded transition
rates
\cite{GuoIEEE:2003,Guo:2006TOP,GuoRider:2006,Guo:2007,Guo:2009,GuoHuangSong:2012,ABPZY:20102,PrietoOHL:2012book,Zhu:2008}.
In comparison, the existence of a Lyapunov function is not needed
in the present paper; Condition \ref{Zy13WeakCon}(c) allows
essentially arbitrarily unbounded transition rates, and thus the
underlying process to be explosive. Part (iii) of Condition
\ref{Zy13GuoYeCon}(c) states the strong continuity of $q(dy|x,a);$
accordingly, the lower semi-continuity of the cost rate $c(x,a)$
is only required in $a\in A(x)$, but the multifunction $A(\cdot)$
needs be compact-valued, which is not required in the present
paper. Finally, one notes that Condition \ref{Zy13GuoYeCon}
implies Condition \ref{Zy13G}.

\section{Proof of Theorem \ref{Zy13MainThm}}\label{Zy13Sec4}
In this section, before proving Theorem \ref{Zy13MainThm}, we
firstly present some auxiliary statements.

Under each Markov policy $\pi\in \Pi_M,$ the process $\xi_t$ is a
Markov jump process \cite{FeinbergS:2013}, and there exists a
transition (sub-probability, in general) function
$p^\pi(u,x,t,dy)$ such that $P^\pi(\xi_t\in
dy|\xi_u)=p^\pi(u,\xi_u,t,dy)$ with $t\ge u\ge 0$ almost surely
with respect to $P^\pi$ \cite{Kuznetsov:1981}. So we formally
define for each $x\in S,$ $u\le t$ and Markov policy $\pi\in
\Pi_M$
\begin{eqnarray}\label{Zy13eqn1}
W^\pi(u,x,t) &:=&\int_u^t \int_S \int_A
c(y,a)\pi(da|y,s)p^\pi(u,x,s,dy)ds
\end{eqnarray}

The next result is a generalization of Theorem 3.4 in
\cite{GuoYe:2010AAP}, which was proved for deterministic
stationary policies only and additionally under Condition
\ref{Zy13GuoYeCon}(a). Since Condition \ref{Zy13GuoYeCon}(a) is
not required in the present article, to be self-contained and for
its potential independent interest, we include this result here,
and present its complete proof in the appendix.

\begin{lemma}\label{Zy13GYLem}

\par\noindent(a) Let a Markov policy $\pi\in \Pi_M$ be fixed. Then the function $W^\pi(u,x,t)$ is the minimal nonnegative
measurable solution to the following inequality
\begin{eqnarray}\label{Zy13Ine1}
v(u,x,t)&\ge& \int_u^t \int_A c(x,a)\pi(da|x,\theta)d\theta~
e^{-\int_u^t \int_A q_x(a)\pi(da|x,\theta)d\theta}\nonumber\\
&&+\int_u^t e^{-\int_u^s \int_A
q_x(a)\pi(da|x,\theta)d\theta}\left\{\int_A
q_x(a)\pi(da|x,s)\int_u^s \int_A
c(x,a)\pi(da|x,\theta)d\theta\right.\nonumber\\
&&\left.+\int_{S\setminus\{x\}} \int_A
q(dy|x,a)\pi(da|x,s)v(s,y,t)\right\}ds.
\end{eqnarray}
\par\noindent(b)
Let a stationary policy $\pi$ be fixed, and suppose there exist a
constant $g\in[0,\infty]$ and a nonnegative measurable function
$h$ on $S$ satisfying the following inequality
\begin{eqnarray*}
g+h(x)\int_A q_x(a)\pi(da|x)\ge \int_A c(x,a)\pi(da|x)+\int_{S\setminus\{x\}}h(y)\int_A q(dy|x,a)\pi(da|x)
\end{eqnarray*}
for each $x\in S.$ Then $g\ge W(x,\pi)$ for each $x\in S$ such that $h(x)<\infty.$
\end{lemma}
\par\noindent\textit{Proof}. See the appendix. $\hfill\Box$

The next lemma, to be used in the proof of Theorem
\ref{Zy13MainThm} below, extends some known results for discounted
CTMDPs in the literature \cite{FeinbergBook:2012,Guo:2007} to
weaker conditions.
\begin{lemma}\label{Zy13discountedLem}
Suppose Condition \ref{Zy13WeakCon}(c) is satisfied. For each
$\alpha>0$, $W_\alpha$ is the minimal nonnegative lower
semi-analytic solution to the equation
\begin{eqnarray}\label{Zy13E1}
v(x)=\inf_{a\in
A(x)}\left\{\frac{c(x,a)}{\alpha+w(x)}+\frac{w(x)}{w(x)+\alpha}\int_S
v(y) \left(\frac{q(dy|x,a)}{w(x)}+I\{x\in dy\}\right) \right\}.
\end{eqnarray}
If additionally Condition \ref{Zy13WeakCon}(a,b) also holds, then
$W_\alpha$ is lower semi-continuous on $S$, and there exists a
deterministic stationary for the discounted CTMDP problem.
\end{lemma}
\par\noindent\textit{Proof}.
See the appendix. $\hfill\Box$

\par\noindent\textbf{Proof of Theorem \ref{Zy13MainThm}.}
Note that under Condition \ref{Zy13G}, $m_\alpha<\infty$ by
Proposition A.5 of \cite{Guo:2009} for all sufficiently small
$\alpha>0$, say, to be specific, for all $0<\alpha\le
\alpha_0<\infty.$ Indeed, Condition \ref{Zy13G} asserts the
existence of some $z\in S$ and policy $\pi\in \Pi$ such that
$W(z,\pi)=\limsup_{t\uparrow\infty} \frac{1}{t}E_z^\pi\left[
\int_0^t \int_A c(\xi_s,a)\pi(da|\omega,s)ds\right]<\infty.$ Thus,
for all sufficiently large $t>0,$ $\int_0^t E_z^\pi\left[\int_A
c(\xi_s,a)\pi(da|\omega,s)\right]ds= E_z^\pi\left[ \int_0^t \int_A
c(\xi_s,a)\pi(da|\omega,s)ds\right]<\infty$. Due to the
nonnegativity of the cost rate $c$, this implies
$E_z^\pi\left[\int_A c(\xi_t,a)\pi(da|\omega,t)\right]<\infty$ for
$t>0$ almost everywhere. Thus the condition of Proposition A.5 in
\cite{Guo:2009} is verified, and we infer from it for that
\begin{eqnarray}\label{ZY13GuoProp}
\limsup_{0<\alpha\downarrow 0}\alpha W_\alpha(z,\pi) \le
W(z,\pi)<\infty,
\end{eqnarray}
and consequently, there exists some $0<\alpha_0<\infty$ such that
$m_\alpha\le W_\alpha(z,\pi)<\infty$ for all $0<\alpha\le
\alpha_0$ as required.

Let
$
g_\alpha:=\alpha m_\alpha.
$
For each $0<\alpha\le \alpha_0,$ we write
$W_\alpha(x)=h_\alpha(x)+m_\alpha$ in (\ref{Zy13E1}) with
$W_\alpha$ in lieu of $v$, and obtain
\begin{eqnarray}\label{Zy13Proof1}
h_\alpha(x)+m_\alpha &=& \inf_{a\in
A(x)}\left\{\frac{c(x,a)}{\alpha+w(x)}+\frac{w(x)}{w(x)+\alpha}\int_S
(h_\alpha(y)+m_\alpha) \left(\frac{q(dy|x,a)}{w(x)}+I\{x\in dy\}\right) \right\}\nonumber\\
&=& \inf_{a\in
A(x)}\left\{\frac{c(x,a)}{\alpha+w(x)}+\frac{w(x)}{w(x)+\alpha}\int_S
h_\alpha(y) \left(\frac{q(dy|x,a)}{w(x)}+I\{x\in
dy\}\right)+\frac{w(x)m_\alpha}{\alpha+w(x)} \right\}.
\end{eqnarray}
It follows from (\ref{Zy13Proof1}) that
\begin{eqnarray}\label{Zy13Proof2}
(w(x)+\alpha)h_\alpha(x)+g_\alpha&=&\inf_{a\in
A(x)}\left\{c(x,a)+w(x)\int_S
h_\alpha(y)\left(\frac{q(dy|x,a)}{w(x)}+I\{x\in
dy\}\right)\right\}.
\end{eqnarray}
Define now
\begin{eqnarray}\label{Zy13Newlyadded}
g:=\limsup_{0<\alpha\downarrow 0} g_\alpha\ge 0,
\end{eqnarray}
which is finite because of (\ref{ZY13GuoProp}), and
\begin{eqnarray}\label{Zy2013Newlyadded2}
h(x):=\liminf_{0<\alpha\downarrow 0,~y\rightarrow x}h_\alpha(y),
\end{eqnarray}
which is finite under Condition \ref{Zy13ConRelative}. It is known
that for each convergent sequence $0<\alpha_n\downarrow 0$ as
$n\rightarrow \infty,$
\begin{eqnarray}\label{Zy13FC}
\sup_{\alpha\in (0,\infty)}\underline{h}_\alpha(x)=h(x)=\liminf_{n\rightarrow\infty, y\rightarrow
x}\underline{h}_{\alpha_n}(y),
\end{eqnarray}
where $ \underline{h}_{\alpha_n}(x):=\liminf_{y\rightarrow x}
H_{\alpha_n}(y) $ with $
H_{\alpha_n}(y):=\inf_{\alpha\in(0,\alpha_n]}h_\alpha(y), $ see
(24) of \cite{Feinberg:2012b} for the first equality in
(\ref{Zy13FC}) and Corollary 1 of \cite{Feinberg:2012b} for the
other. The above three functions are all measurable; in fact, the
functions $\underline{h}_\alpha$ and $h$ are lower semi-continuous
on $S$, see Lemma 5.13.4 of \cite{Berberian:1999} and Lemma 4.2 of
\cite{CavazosCadena:2010}, respectively. Note that by their
definitions
\begin{eqnarray}\label{Zy13Proof3}
h_\beta(x)\ge H_\beta(x)\ge H_\alpha(x)\ge \underline{h}_\alpha(x)
\end{eqnarray}
for each $x\in S$ and $\alpha\ge \beta>0.$

Let $\epsilon>0$ be arbitrarily fixed. Then by the definition of
the constant $g$ (see (\ref{Zy13Newlyadded})), there exists
$0<\alpha_1\le \alpha_0$ such that for each
$\alpha\in(0,\alpha_1],$ $g\ge g_\alpha-\epsilon.$ It follows from
this and (\ref{Zy13Proof2}) that for each $0<\alpha\le \alpha_1,$
\begin{eqnarray*}
(w(x)+\alpha)h_\alpha(x)+g+\epsilon\ge\inf_{a\in
A(x)}\left\{c(x,a)+w(x)\int_S
h_\alpha(y)\left(\frac{q(dy|x,a)}{w(x)}+I\{x\in
dy\}\right)\right\},
\end{eqnarray*}
which, together with
(\ref{Zy13Proof3}), leads to that for each $0<\beta\le \alpha\le
\alpha_1,$
\begin{eqnarray*}
(w(x)+\beta)h_\beta(x)+g+\epsilon\ge\inf_{a\in
A(x)}\left\{c(x,a)+w(x)\int_S
H_\alpha(y)\left(\frac{q(dy|x,a)}{w(x)}+I\{x\in
dy\}\right)\right\},
\end{eqnarray*}
and thus by the definition of $H_\alpha$, the above relation and (\ref{Zy13Proof3}) again,
\begin{eqnarray}\label{Zy13proof4}
(w(x)+\alpha)H_\alpha(x)+g+\epsilon&\ge&
\inf_{a\in
A(x)}\left\{c(x,a)+w(x)\int_S
H_\alpha(y)\left(\frac{q(dy|x,a)}{w(x)}+I\{x\in
dy\}\right)\right\}   \nonumber \\
 &\ge&
\inf_{a\in
A(x)}\left\{
c(x,a)+
w(x)\int_S
\underline{h}_\alpha(y)\left(\frac{q(dy|x,a)}{w(x)}+I\{x\in
dy\}\right)\right\}
\end{eqnarray}
for each $0<\alpha\le \alpha_1.$ Under Condition
\ref{Zy13WeakCon}, the stochastic kernel
$\frac{q(dy|x,a)}{w(x)}+I\{x\in dy\}$ is weakly continuous, which,
together with the lower semi-continuity of $\underline{h}_\alpha$
(as explained earlier), implies that
\begin{eqnarray*}
w(x)\int_S
\underline{h}_\alpha(y)\left(\frac{q(dy|x,a)}{w(x)}+I\{x\in
dy\}\right)
\end{eqnarray*}
defines a lower semi-continuous function on $S$. As a result,
$c(x,a)+ w(x)\int_S
\underline{h}_\alpha(y)\left(\frac{q(dy|x,a)}{w(x)}+I\{x\in
dy\}\right)$ is $\mathbb{K}$-inf-compact because so is the cost
rate $c$ and that $\underline{h}_\alpha(x)\ge 0$ for each $x\in
S,$ see Lemma \ref{Zy13NEW} in the appendix. Therefore, one can
infer from Lemma \ref{Zy13BergeF} in the appendix for the lower
semi-continuity on $S$ of the expression in the second line of
(\ref{Zy13proof4}). Following from this and upon taking the
corresponding lower limit on the both sides of (\ref{Zy13proof4}),
one obtains
\begin{eqnarray*}
(w(x)+\alpha)\underline{h}_\alpha(x)+g+\epsilon\ge
\inf_{a\in
A(x)}\left\{
c(x,a)+
w(x)\int_S
\underline{h}_\alpha(y)\left(\frac{q(dy|x,a)}{w(x)}+I\{x\in
dy\}\right)\right\}
\end{eqnarray*}
for each $0<\alpha\le \alpha_1.$  Now the first equality of
(\ref{Zy13FC}) and the above inequality imply
\begin{eqnarray}\label{Zy13Proof5}
(w(x)+\alpha) h(x)+g+\epsilon\ge \inf_{a\in A(x)}\left\{ c(x,a)+
w(x)\int_S
\underline{h}_\alpha(y)\left(\frac{q(dy|x,a)}{w(x)}+I\{x\in
dy\}\right)\right\}
\end{eqnarray}
for each $0<\alpha\le \alpha_1.$ By the
$\mathbb{K}$-inf-compactness of the expression inside the
parenthesis on the right side of (\ref{Zy13Proof5}) (as explained
earlier) and Lemma \ref{Zy13BergeF}, for each $0<\alpha\le
\alpha_1,$ there exists some $a_\alpha\in A(x)$ such that
\begin{eqnarray}\label{Zy13Proof6}
 (w(x)+\alpha) h(x)+g+\epsilon&\ge&
\inf_{a\in A(x)}\left\{
c(x,a)+
w(x)\int_S
\underline{h}_\alpha(y)\left(\frac{q(dy|x,a)}{w(x)}+I\{x\in
dy\}\right)\right\}  \nonumber\\
&=&
c(x,a_\alpha)+
w(x)\int_S
\underline{h}_\alpha(y)\left(\frac{q(dy|x,a_\alpha)}{w(x)}+I\{x\in
dy\}\right).
\end{eqnarray}
Now let $x\in S$ be arbitrarily fixed, and take $\alpha_1\ge
\alpha_n\downarrow 0$ ($\alpha_n>0$). Under Condition
\ref{Zy13ConRelative} the expression on the left side of
inequality (\ref{Zy13Proof6}) is finite (recall
(\ref{Zy2013Newlyadded2}) for the definition of the function $h$).
Considering (\ref{Zy13Proof6}) with $\alpha_n$ replacing $\alpha$
therein, it follows from the definition of the
$\mathbb{K}$-inf-compactness that the sequence $\{a_{\alpha_n}\}$
admits a limit point $a^\ast\in A(x).$ Taking the lower limit on
the both sides of (\ref{Zy13Proof6}) along the specified sequence
$\alpha_1\ge \alpha_n\downarrow 0$ ($\alpha_n>0$), we see
 \begin{eqnarray}\label{Zy13Proof7}
 w(x)h(x)+g+\epsilon&\ge&
 c(x,a^\ast)+
 w(x)\liminf_{n\rightarrow\infty}\int_S
 \underline{h}_{\alpha_n}(y)\left(\frac{q(dy|x,a_{\alpha_n})}{w(x)}+I\{x\in
 dy\}\right)\nonumber\\
&\ge& c(x,a^\ast)+ w(x)\int_S
h(y)\left(\frac{q(dy|x,a^\ast)}{w(x)}+I\{x\in
dy\}\right)\nonumber\\
&\ge&
\inf_{a\in A(x)}\left\{c(x,a)+
w(x)\int_S
h(y)
\left(\frac{q(dy|x,a)}{w(x)}+I\{x\in
dy\}\right)\right\},
 \end{eqnarray}
where for the first inequality the finiteness of $h(x)$ and the
lower semi-continuity of the term inside the parenthesis on the
right side of (\ref{Zy13Proof5}) are used; and the second
inequality follows from (\ref{Zy13FC}), the weak continuity of the
underlying stochastic kernel, and the generalized Fatou's lemma,
see Lemma \ref{Zy13FatouF} in the appendix or Lemma 4.2 of
\cite{CavazosCadena:2010}. That the inequality in
(\ref{Zy13OptIne}) is satisfied by the constant $g$ and the
nonnegative real-valued lower semi-continuous function $h$ follows
from (\ref{Zy13Proof7}) and the arbitrariness of $\epsilon>0$.
Regarding the existence of a measurable selector $\varphi^\ast$
satisfying the equality in (\ref{Zy13OptIne}), one can refer to
Lemma \ref{Zy13BergeF}; recall that the term in the parenthesis in
(\ref{Zy13OptIne}) is $\mathbb{K}$-inf-compact. We prove the rest
of this statement as follows. Let $\varphi^\ast$ be any measurable
selector satisfying the equality in (\ref{Zy13OptIne}). By the
finiteness of $h(x)$, (\ref{Zy13OptIne}) and Lemma
\ref{Zy13GYLem},
\begin{eqnarray}\label{Zy13Proof8}
g\ge W(x,\varphi^\ast) \ge \inf_{\pi\in \Pi}W(x,\pi).
\end{eqnarray}
For the opposite direction, let $x\in S$ be arbitrarily fixed.
Since $g<\infty,$ we see $\inf_{\pi\in \Pi}W(x,\pi)<\infty.$ Fix
arbitrarily some (possibly $x$-dependent) policy $\pi$ such that
$W(x,\pi)<\infty.$ Now as in the argument for (\ref{ZY13GuoProp})
with $z$ being replaced by $x$ in the beginning of this proof, we
see $\limsup_{0<\alpha\downarrow 0}\alpha W_\alpha(x)\le
W(x,\pi)<\infty,$ which together with the arbitrariness of the
policy $\pi$ and the fact that $g=\limsup_{0<\alpha\downarrow
0}\alpha \inf_{x\in S}W_\alpha(x)\le \limsup_{0<\alpha\downarrow
0}\alpha W_\alpha(x)$ (recalling here the definition of $g$ given
by (\ref{Zy13Newlyadded})), leads to $\inf_{\pi\in \Pi}
W(x,\pi)\ge g.$ Thus, we see the validity of (\ref{Zy13Proof8})
with inequalities being replaced by equalities. It follows from
the arbitrariness of $x\in S$ that the policy $\varphi^\ast$ is
optimal. The proof is now completed. $\hfill\Box$

\section{Conclusion}\label{Zy13Sec5}
To sum up, for a CTMDP in Borel state and action spaces with a
nonnegative cost rate, the existence of a deterministic stationary
average optimal policy is proved with weakly continuous transition
rates. Our conditions allow the controlled process to be explosive
(i.e., the transition rates are essentially arbitrarily
unbounded). In addition, following the neat generalization of the
Berge theorem \cite{FeinbergKZ:2013}, the condition on the
admissible action spaces has been further relaxed as compared with
the previous literature.

\section*{Appendix}

\begin{definition}\label{Zy13DefLSA}
The collection of analytic subsets of a nonempty Borel space $S$
is the collection of images of measurable subsets of $Y$ under all
measurable mappings from $Y$ into $S$, where $Y$ is an uncountable
Borel space. A function $f$ on the nonempty Borel space $S$ is
called lower semi-analytic if for each $\epsilon
\in(-\infty,\infty)$, the set $\{x\in S: f(x)<\epsilon\}$ is
analytic. A function $f$ is called upper semi-analytic if $-f$ is
lower semi-analytic.
\end{definition}
See more details about the above definition in Chapter 7 of
\cite{Bertsekas:1978}.

The next lemma comes from \cite{FeinbergKZ:2013}, see Theorems 1.2
and 3.3 therein, where the more general statements are
established.
\begin{lemma}\label{Zy13BergeF}
Suppose a function $g$ on the nonempty Borel space $K=\{(x,a):x\in
S, a\in A(x)\}$ is $\mathbb{K}$-inf-compact. Then $\inf_{a\in
A(x)} g(x,a)$ defines a lower semi-continuous function in $x\in
S$. Furthermore, there is a measurable mapping $\varphi^\ast$ from
$S$ to $A,$ whose graph is contained in $K,$ such that $\inf_{a\in
A(x)} g(x,a)=g(x,\varphi^\ast(x))$ for each $x\in S.$
\end{lemma}

The following lemma summarizes some facts about
$\mathbb{K}$-inf-compact functions, which are used frequently in
the proofs in this paper.
\begin{lemma}\label{Zy13NEW}
Let $c$ be a $\mathbb{K}$-inf-compact function on $K.$ If $v$ is a
nonnegative lower semi-continuous function on $K$, then $c+v$ is
also $\mathbb{K}$-inf-compact on $K.$ If $u$ is a continuous
real-valued function on $S$ such that $u(x)> 0$ for each $x\in S,$
then $\frac{c(x,a)}{u(x)}$ defines a $\mathbb{K}$-inf-compact
function on $K.$
\end{lemma}
\par\noindent\textit{Proof.} We only verify the second part. Clearly $\frac{c(x,a)}{u(x)}$
is lower semi-continuous on $K$. Now suppose $S\ni x_n\rightarrow
x\in S$ and $a_n\in A(x_n)$ such that there is some real constant
$M>0$ such that $ \frac{c(x_n,a_n)}{u(x_n)}\le M,$ i.e.,
$c(x_n,a_n)\le M u(x_n)$ for all $n.$ Since $u$ is continuous and
the set $X:=\bigcup_{n=0}^\infty \{x_n\}\bigcup\{x\}$ is compact
in $S$, we further infer from the previous inequality for that $
c(x_n,a_n)\le M \sup_{y\in X}u(y)<\infty$ for all $n.$ Now it
follows from the $\mathbb{K}$-inf-compactness of the function $c$
that there exists a limit point $a\in A(x)$ for the sequence
$\{a_n\}$, as required. $\hfill\Box$

The following statement is known as the generalized Fatou's lemma
\cite{CavazosCadena:2010,CostaD:2012,FeinbergKZ:2013b}. A detailed
proof with more general statements is available at
\cite{FeinbergKZ:2013b}.
\begin{lemma}\label{Zy13FatouF}
Suppose a sequence of probability measures $Q_n$ on the nonempty
Borel space ${\cal B}(S)$ is weakly convergent to the probability
measure $Q$ on ${\cal B}(S).$ Then for each sequence of
nonnegative functions $g_n$ on $S$, it holds that
$\int_{S}(\liminf_{n\rightarrow\infty,x\rightarrow y}
g_n(x))Q(dy)\le \liminf_{n\rightarrow\infty}\int_S g_n(y)Q_n(dy).$
\end{lemma}

\par\noindent\textbf{Proof of Lemma \ref{Zy13GYLem}.}
(a) For simplicity, throughout the proof of this lemma, we omit
the fixed policy $\pi$ from indications, and introduce the
following
 notations
 \begin{eqnarray*}
 c(x,s):=\int_A c(x,a)\pi(da|x,s),~q_x(s):=\int_A
 q_x(a)\pi(da|x,s),~q(dy|x,s):=\int_A q(dy|x,a)\pi(da|x,s).
 \end{eqnarray*}
Furthermore, if $c(x,s)$, $q_x(s)$ and $q(dy|x,s)$ in the above are $s$-independent, as in the case of a stationary policy, we
omit $s$ from the arguments.

It is known \cite{FeinbergS:2013} that the transition function
$p(u,x,t,dy)$ can be constructed iteratively by
$
\sum_{k=0}^n p_k(u,x,t,dy)\uparrow p(u,x,t,dy)
$
as $n\uparrow\infty,$ where the convergence is set-wise, and for
each $\Gamma\in {\cal B}(S)$
\begin{eqnarray*}
p_0(u,x,t,\Gamma)&:=&I\{x\in \Gamma\}e^{-\int_u^t q_x(s)ds};\\
p_k(u,x,t,\Gamma)&:=&\int_u^t \int_{S\setminus\{x\}} e^{-\int_u^s
q_x(\theta)d\theta } q(dy|x,s) p_{k-1}(s,y,t,\Gamma)ds.
\end{eqnarray*}
It follows from this, the nonnegativity of the cost rate $c$ and
the monotone convergence theorem, see Theorem 2.1 in
\cite{Hernandez-Lerma:2000}, that
$
m_n(u,x,t):=\int_u^t \int_S c(y,s) \sum_{k=0}^n
p_n(u,x,s,dy)ds\uparrow W(u,x,t)
$
as $n\uparrow\infty$, see (\ref{Zy13eqn1}).

We verify firstly that $W(u,x,t)$ satisfies (\ref{Zy13Ine1}) with
equality as follows. By the iterative definitions of the
transition functions $p_n$,
\begin{eqnarray*}
m_n(u,x,t)&:=&\int_u^t \int_S c(y,s)\sum_{k=0}^n p_k(u,x,s,dy)
ds\nonumber\\
&=&m_0(u,x,t)+\int_u^t \int_S c(y,s)\sum_{k=1}^n p_k(u,x,s,dy)ds\nonumber\\
&=&m_0(u,x,t)+\int_u^t \int_S c(y,s)\sum_{k=1}^n \int_u^s
\int_{S\setminus \{x\}} e^{-\int_u^r q_x(\theta)d\theta}
q(dz|x,r)p_{k-1}(r,z,s,dy) dr~ds  \nonumber  \\
&=&m_0(u,x,t)+\int_u^t \int_S c(y,s)\sum_{k-1=0}^{n-1} \int_r^t   \int_{S\setminus \{x\}} e^{-\int_u^r q_x(\theta)d\theta}   q(dz|x,r)p_{k-1}(r,z,s,dy) ds~dr\nonumber\\
&=&m_0(u,x,t)+\int_u^t e^{-\int_u^r q_x(\theta)d\theta} \int_{S\setminus\{x\}} q(dz|x,r)m_{n-1}(r,z,t)dr,
\end{eqnarray*}
where the last two inequalities follow from the legal interchange
of the order of integrations. Integration by parts gives $
m_0(u,x,t)=e^{-\int_u^t q_x(\theta)d\theta}\int_u^t
c(x,\theta)d\theta+\int_u^t \int_u^s c(x,\theta)d\theta
e^{-\int_u^s q_x(\theta)d\theta}q_x(s)ds. $ It thus follows that
\begin{eqnarray}\label{Zy13eq2}
m_n(u,x,t)&=&\int_u^t c(x,\theta)d\theta e^{-\int_u^t q_x(\theta)d\theta} \nonumber\\
&&+\int_u^t e^{-\int_u^s q_x(\theta)d\theta}\left\{q_x(s)\int_u^s c(x,\theta)d\theta+\int_{S\setminus \{x\}} q(dy|x,s)m_{n-1}(s,y,t) \right\}ds.
\end{eqnarray}
By the standard monotone convergence theorem, passing to the limit as $n\uparrow\infty$ on the both sides of the above equality gives
\begin{eqnarray*}
W(u,x,t)&=&\int_u^t c(x,\theta)d\theta e^{-\int_u^t q_x(\theta)d\theta} \nonumber\\
&&+\int_u^t e^{-\int_u^s q_x(\theta)d\theta}\left\{q_x(s)\int_u^s c(x,\theta)d\theta+\int_{S\setminus \{x\}} q(dy|x,s)W(s,y,t) \right\}ds.
\end{eqnarray*}

For the minimality of $W(u,x,t)$ as a nonnegative measurable
solution to inequality (\ref{Zy13Ine1}), suppose that there is
another nonnegative measurable solution $v(u,x,t)$ to inequality
(\ref{Zy13Ine1}). Thus, $v(u,x,t)\ge m_0(u,x,t)$. Now an inductive
argument based on (\ref{Zy13eq2}) and the fact that $v$ satisfies
(\ref{Zy13Ine1}) implies $v(u,x,t)\ge m_n(u,x,t)$ for each
$n=0,1,\dots,$ which, together with the fact that $m_n\uparrow W$
point-wise as $n\uparrow\infty,$ leads to that $v(u,x,t)\ge
W(u,x,t)$ as desired.

(b) Suppose a stationary policy $\pi$ is fixed, and there exist a
constant $g$ and a nonnegative measurable function $h$ on $S$ as
in the statement. Without loss of generality, we assume that
$g<\infty$ for otherwise the statement holds automatically. It is
well known, or otherwise follows from the construction of the
transition function $p(u,x,t,dy)$ above that under the stationary
policy, $p(u,x,t,dy)$ depends on $u$ and $t$ only through the time
increment $t-u$, and the underlying Markov jump process $\xi_t$ is
homogeneous, and thus $W(u,x,t)=E_x\left[\int_0^{t-u} c(\xi_ s)ds
 \right]=:\tilde{W}(x,t-u)$, see Theorem 2.2 of \cite{FeinbergS:2013}; recall the agreement that
the (stationary) policy $\pi$ is omitted from indication in this
proof. It follows from this and part (a) specialized to a
stationary policy and $u=0$, that $\tilde{W}(x,t)$ is the minimal
nonnegative measurable solution to the inequality
\begin{eqnarray*}
\tilde{W}(x,t)&\ge&  c(x)t e^{- q_x t}+\int_0^t e^{- q_x
s}\left\{q_x  ~ c(x)s  +\int_{S\setminus \{x\}} q(dy|x )\tilde{W}(
y,t-s) \right\}ds.
\end{eqnarray*}
Now it can be verified, based on the definitions of the constant $g$ and the function $h$, that
the above inequality is satisfied with $h(x)+gt$ in lieu of $\tilde{W}(x,t)$. Consequently,
$
h(x)+gt \ge \tilde{W}(x,t)
$
by part (a) of this lemma.
At $x\in S$ such that $h(x)<\infty,$ dividing the both sides of the previous inequality and then passing to the upper limit as $t\rightarrow\infty$ yields the statement.
$\hfill\Box$
\bigskip

\par\noindent\textbf{Proof of Lemma \ref{Zy13discountedLem}.} Let $\alpha>0$ be arbitrarily fixed. It
is known that the value function $W_\alpha$ for the discounted
CTMDP problem is the minimal nonnegative lower semi-analytic
solution to the equation
\begin{eqnarray}\label{Zy13E2}
v(x)=\inf_{a\in
A(x)}\left\{\frac{c(x,a)}{\alpha+q_x(a)}+\int_{S\setminus
\{x\}}v(y) \frac{q(dy|x,a)}{\alpha+q_x(a)}\right\}=:\tilde{T}\circ
v(x),
\end{eqnarray}
see Theorem 5.5.5 in \cite{FeinbergBook:2012}. For the first part
of this lemma, it remains to recognize that the two equations
(\ref{Zy13E1}) and (\ref{Zy13E2}) admit the same minimal
nonnegative solution. Below, in spite that the argument is
trivial, we briefly verify this relation because first, a similar
relation between equation (\ref{Zy13E1}) and another equation
similar to (\ref{Zy13E2}) was falsely claimed without proofs in
\cite{PiunovskiyZY:2012}, see equation (8) therein, and second, it
is easy to construct examples to show that equations
(\ref{Zy13E1}) and (\ref{Zy13E2}) are not equivalent; indeed,
there can be solutions to (\ref{Zy13E1}), which do not satisfy
(\ref{Zy13E2}). For brevity, we write (\ref{Zy13E1}) as $v=T\circ
v$ with
$
T\circ v(x):=\inf_{a\in
A(x)}\left\{\frac{c(x,a)}{\alpha+w(x)}+\frac{w(x)}{w(x)+\alpha}\int_S
v(y) \left(\frac{q(dy|x,a)}{w(x)}+I\{x\in dy\}\right)\right\}.
$
Firstly, consider the minimal nonnegative solution $u$ to
(\ref{Zy13E2}), and let $x\in S$ be arbitrarily fixed. If
$u(x)=\infty,$ then $T\circ u(x)=\infty=u(x)$ (recalling the
convention of $\infty-\infty:=\infty$). Now suppose $u(x)<\infty.$
Then it follows that $u(x)\le
\frac{c(x,a)}{\alpha+w(x)}+\frac{w(x)}{w(x)+\alpha}\int_S u(y)
\left(\frac{q(dy|x,a)}{w(x)}+I\{x\in dy\}\right)$ for each $a\in
A(x).$  Let $\delta>0$ be arbitrarily fixed, and take any
$0<\epsilon<\delta.$ Then there exists some $a_\delta\in A(x)$
such that $u(x)+\epsilon\ge
\frac{c(x,a_\delta)}{\alpha+q_x(a_\delta)}+\int_{S\setminus
\{x\}}u(y) \frac{q(dy|x,a_\delta)}{\alpha+q_x(a_\delta)}$ so that
$u(x)+\delta
>u(x)+\frac{\epsilon (\alpha+q_x(a_\delta))}{\alpha+w(x)}\ge
\frac{c(x,a_\delta)}{\alpha+w(x)}+\frac{w(x)}{w(x)+\alpha}\int_S
u(y) \left(\frac{q(dy|x,a_\delta)}{w(x)}+I\{x\in dy\}\right).$
Since $\delta>0$ is arbitrarily fixed, we see that $u(x)=T\circ
u(x)$. Thus, $u\ge v$ with $v$ being the minimal nonnegative
solution to (\ref{Zy13E1}). For the opposite direction, note that
if $v(x)=\infty,$ then $v(x)\ge \tilde{T}\circ v(x)$. Suppose now
$v(x)<\infty.$ Then for each $a\in A(x)$, $v(x)\le
\frac{c(x,a)}{\alpha+w(x)}+\frac{w(x)}{w(x)+\alpha}\int_S v(y)
\left(\frac{q(dy|x,a)}{w(x)}+I\{x\in dy\}\right),$ and so $v(x)\le
\frac{c(x,a)}{\alpha+q_x(a)}+\int_{S\setminus\{x\}}
v(y)\frac{q(dy|x,a)}{q_x(a)+\alpha}.$ Let $\delta>0$ be
arbitrarily fixed, and choose $\epsilon>0$ such that
$\frac{\epsilon(\alpha+w(x))}{\alpha}<\delta.$ Since $v$ satisfies
(\ref{Zy13E1}), there exists some $a_\delta\in A(x)$ such that $
v(x)\ge
\frac{c(x,a_\delta)}{\alpha+w(x)}+\frac{1}{\alpha+w(x)}\int_S
v(y)q(dy|x,a_\delta) +\frac{w(x)v(x)}{\alpha+w(x)}-\epsilon. $
Simple rearrangements of this inequality further lead to $v(x)\ge
\frac{c(x,a_\delta)}{\alpha+q_x(a_\delta)}+\frac{1}{\alpha+q_x(a_\delta)}\int_{S\setminus\{x\}}
v(y)q(dy|x,a_\delta)-\delta.$ Thus, $v(x)\ge \tilde{T}\circ v(x)$.
It follows from this and Proposition 9.10 of \cite{Bertsekas:1978}
that $u\le v,$ and thus $u=v$ (recalling the opposite direction of
the previous inequality being established earlier). The first part
of this lemma is proved.

Next, we observe that according to the first part of this lemma
and Proposition 9.16 of \cite{Bertsekas:1978}, $W_\alpha$ is also
given by the value function of a DTMDP with the total undiscounted
cost criterion specified by the following primitives. The state
space is $S\bigcup \{x_\infty\}$; the action space is $A\bigcup
\{a_\infty\}$; the admissible action space is $A(x)$ for each
$x\in S$ with $A(x_\infty)=\{a_\infty\};$ the transition
probability is given by $Q(\Gamma|x,a):=
\frac{w(x)}{w(x)+\alpha}\left(\frac{q(\Gamma|x,a)}{w(x)}+I\{x\in
\Gamma\}\right)$ for each $x\in S$, $a\in A(x)$ and $\Gamma\in
{\cal B}(S),$
$Q(\{x_\infty\}|x,a):=I\{x=x_\infty,a=a_\infty\}+I\{x\in S,a\in
A(x)\}(1-Q(S|x,a));$ and finally, the cost function is $I\{x\in
S,a\in A(x)\}\frac{c(x,a)}{\alpha+w(x)}.$ Here we recall that
$x_\infty\notin S$ and $a_\infty\notin A$ are two isolated points.
Under Condition \ref{Zy13WeakCon}, one can verify that the
transition probability $Q(dy|x,a)$ is weakly continuous, i.e., for
each bounded continuous function $f$ on $S\bigcup \{x_\infty\}$,
$\int_{S\bigcup\{x_\infty\}} f(y)Q(dy|x,a)$ is continuous in $x\in
S\bigcup\{x_\infty\}$ and $a\in A(x);$ and the cost function is
$\mathbb{K}$-inf-compact, see Lemma \ref{Zy13NEW}. Denote the
value function for this DTMDP problem with the total undiscounted
cost criterion also by $W_\alpha$. Below, to be self-contained, we
verify that $W_\alpha$ can be constructed using the value
iteration algorithm under Condition \ref{Zy13WeakCon}. Let
$v_0(x):=0$ and $v_n(x):= \inf_{a\in
A(x)}\left\{\frac{c(x,a)}{\alpha+w(x)} +
\frac{w(x)}{w(x)+\alpha}\int_S
v_{n-1}(y)\left(\frac{q(dy|x,a)}{w(x)}+I\{x\in
dy\}\right)\right\}$ for each $x\in S$, whereas $v_n(x_\infty):=0$
for each $n=0,1,2,\dots.$ Under Condition \ref{Zy13WeakCon}, since
the transition probability is weakly continuous and the cost
function is $\mathbb{K}$-inf-compact, by Lemma \ref{Zy13BergeF},
$v_n$ is lower semi-continuous for each $n=0,1,\dots.$
Furthermore, the sequence $\{v_n\}$ is increasing, so that we
formally define $v_\infty(x):=\lim_{n\uparrow\infty}v_n(x)$, which
is thus also lower semi-continuous. Let $x\in S$ be arbitrarily
fixed. It is easy to see from the monotone convergence theorem
that $v_\infty(x)\le \frac{c(x,a)}{\alpha+w(x)}+\int_S v_\infty(y)
Q(dy|x,a)$ for each $a\in A(x)$, and thus $v_\infty(x)\le
\inf_{a\in A(x)}\left\{\frac{c(x,a)}{\alpha+w(x)}+\int_S
v_\infty(y) Q(dy|x,a)\right\}.$ For the opposite direction,
without loss of generality, we assume that $v_\infty<\infty.$ For
each fixed $m\le n-1,$
\begin{eqnarray}\label{Zy13Inter}
v_\infty(x)&\ge& v_n(x)= \inf_{a\in
A(x)}\left\{\frac{c(x,a)}{\alpha+w(x)}+\int_S v_{n-1}(y)
Q(dy|x,a)\right\}\nonumber\\
&=& \frac{c(x,a_n)}{\alpha+w(x)}+\int_S
v_{n-1}(y)Q(dy|x,a_n)\ge\frac{c(x,a_n)}{\alpha+w(x)}+\int_S
v_{m}(y)Q(dy|x,a_n),
\end{eqnarray}
where $a_n\in A(x)$ are the corresponding minimizers, whose
existence is ensured by Lemma \ref{Zy13BergeF}, and the last
inequality is due to that $\{v_n\}$ is an increasing sequence.
Having noted that $\frac{c(x,a)}{\alpha+w(x)}+\int_S
v_{m}(y)Q(dy|x,a)$ is $\mathbb{K}$-inf-compact, and
$v_\infty(x)<\infty,$ we see that the sequence $\{a_n\}$ admits
some limit point $a^\ast\in A(x).$ Assume without loss of
generality that $a_n\rightarrow a^\ast$ for otherwise one can take
the corresponding subsequence. By passing to the limit as
$n\rightarrow\infty$ on the both sides of (\ref{Zy13Inter}) and
the lower semi-continuity of the involved functions, we obtain
$v_\infty(x)\ge \frac{c(x,a^\ast)}{\alpha+w(x)}+\int_S
v_{m}(y)Q(dy|x,a^\ast).$ Further passing to the limit as
$m\rightarrow\infty$ on the both sides of the above inequality
yields $v_\infty(x)\ge \frac{c(x,a^\ast)}{\alpha+w(x)}+\int_S
v_{\infty}(y)Q(dy|x,a^\ast)\ge \inf_{a\in A(x)}\left\{ \frac{c(x,a
)}{\alpha+w(x)}+\int_S v_{\infty}(y)Q(dy|x,a )\right\}.$ Hence, in
combination with the other direction as proved earlier, we see
that $v_\infty$ is a nonnegative measurable (in fact, lower
semi-continuous) solution to (\ref{Zy13E1}). This, by virtue of
Proposition 9.16 of \cite{Bertsekas:1978}, shows
$v_\infty(x)=W_\alpha(x),$ and thus the lower semi-continuity of
$W_\alpha$ follows. Consequently, there exists a deterministic
stationary policy $\varphi^\ast$ such that
\begin{eqnarray*}
W_\alpha(x)&=&\inf_{a\in
A(x)}\left\{\frac{c(x,a)}{\alpha+w(x)}+\frac{w(x)}{w(x)+\alpha}\int_S
W_\alpha(y) \left(\frac{q(dy|x,a)}{w(x)}+I\{x\in dy\}\right) \right\}\\
&=&\frac{c(x,\varphi^\ast(x))}{\alpha+w(x)}+\frac{w(x)}{w(x)+\alpha}\int_S
W_\alpha(y) \left(\frac{q(dy|x,\varphi^\ast(x))}{w(x)}+I\{x\in dy\}\right) .
\end{eqnarray*}
Evidently, this policy satisfies
$W_\alpha(x)=W_\alpha(x,\varphi^\ast)$. $\hfill\Box$
\bigskip



\par\noindent\textbf{Acknowledgement.} The author is thankful to the helpful
comments and remarks received from the referee and editor.


\begin{thebibliography}{99}
\footnotesize
\bibitem{Berberian:1999} {\sc Berberian, S.} (1999). {\it
Fundamentals of Real Analysis}. Springer, New York.
\bibitem{Bertsekas:1978} {\sc Bertsekas, D. and Shreve, S.} (1978). {\it Stochastic Optimal Control}. Academic Press, New York.
\bibitem{CavazosCadena:1991} {\sc Cavazos-Cadena, R.} (1991). A counter example on the optimality
equation in Markov decision chains with the average cost
criterion. {\it Syst. Control Lett.} {\bfseries 16}, 387-392.
\bibitem{CavazosCadena:2010} {\sc Cavazos-Cadena, R. and Salem-Silva,
F.} (2010). {The discunted method and eqiivalence of average criteria for risk-sensitive Markov decision processes on Borel spaces}, {\it Appl. Math. Optim.} {\bfseries 61}, 167-190.
\bibitem{CostaD:2012} {\sc Costa, O. and Dufour, F.} (2012).
Average control of Markov decision processes with Feller transition probabilities and general acton spaces. {\it J. Math. Anal. Appl.} {\bfseries 396}, 58-69.


\bibitem{FeinbergLewis:2007} {\sc Feinberg,E. and Lewis, M.}
(2007). Optimality inequalities for average cost Markov decisio
processes and the stochastic cash balance problem. {\it Math.
Oper. Res.} {\bfseries 32}, 769-783.
\bibitem{FeinbergBook:2012} {\sc Feinberg, E.} (2012). {Reduction of discounted continuous-time MDPs with unbounded jump and reward rates to discrete-time total-reward
MDPs}, in {\it  Optimization, Control, and Application of
Stochastic Systems}, 77-97, Hern{\'a}ndez-Hern{\'a}ndez, D. and
Minjarez-Sosa, A. (Eds), Birkhauser.


\bibitem{Feinberg:2012b} {\sc Feinberg, E., Kasyanov, P. and Zadoianchuk,
N.} (2013). Average-cost Markov decision processes with weakly
continuous transition probabilities. {\it Math. Oper. Res.}
{\bfseries 37}, 591-607.


\bibitem{FeinbergKZ:2013} {\sc Feinberg, E., Kasyanov, P. and
Zadoianchuk, N.} (2013). Berge's theorem for noncompact image sets. {\it J. Math. Anal. Appl.} {\bfseries 397}, 255-259.

\bibitem{FeinbergKZ:2013b} {\sc Feinberg, E., Kasyanov, P. and
Zadoianchuk, N.} (2013). Fatou's lemma for weakly convergent
probabilities. Preprint, Department of Applied Mathematics and
Statistics, State University of New York at Stony Brook, available
at arxiv:1206.4073v2.

\bibitem{FeinbergS:2013}{\sc Feinberg, E., Mandava, M. and Shiryaev,
A.} (2013). {On solutions of Kolmogorov's equations for jump
Markov processes}. {\it J. Math. Anal. Appl.} {\bfseries 411},
261-270.


\bibitem{Gihman:1975} {\sc Gihman, I. and Skorohod, A.} (1975). {\it The Theory of Stochastic Processes
II}. Springer, Berlin.

\bibitem{GuoLiu:2001} {\sc Guo, X. and Liu, K.} (2001). A note on
optimaluty conditions for continuous-time Markov decision processes with avrage cost criterion. {\it IEEE Trans. Automat. Control} {\bfseries 46}, 1984-1989.
\bibitem{GuoIEEE:2003} {\sc Guo, X. and Hern\'{a}ndez-Lerma, O.} (2003). Drift and monotonicity conditions for continuous-time controlled Markov chains with an average
criterion. {\it IEEE Trans. Automat. Control} {\bfseries 48},
236-245.
\bibitem{Guo:2006TOP} {\sc Guo, X.P., Hern{\'a}ndez-Lerma, O. and Prieto-Rumeau, T.} (2006). {A survey of recent results on continuous-time Markov decision processes}, {Top.} {\bfseries 14}, 177--257.
\bibitem{GuoRider:2006} {\sc Guo, X. and Rieder, U.} (2006).
Average optimality for continuous-time Markov decision processes
in Polish spaces. {\it Ann. Appl. Probab.} {\bfseries 16},
730-756.
\bibitem{Guo:2007} {\sc Guo, X.} (2007). Continuous-time Markov decision processes with discounted rewards: the case of Polish spaces. {\it Math. Oper. Res.} {\bfseries32}, 73-87.

\bibitem{Guo:2009} {\sc Guo, X. and Hern\'{a}ndez-Lerma, O.} (2009). {\em {C}ontinuous-{t}ime {M}arkov {D}ecision {P}rocesses: {T}heory and {A}pplications}. Springer, Heidelberg.


\bibitem{GuoYe:2010AAP} {\sc Guo, X. and Ye, L.} (2010). New discount and average optimality conditions for continuous-time Markov decision processes. {\it Adv. Appl. Probab.} {\bfseries
42}, 953-985.
\bibitem{GuoHuangSong:2012} {\sc Guo, X., Huang, Y. and Song, X.}
(2012). Linear programming and constrained average optimality for
general continuous-time Markov decision processes in
history-dependent policies. {\it SIAM J. Control Optim.}
{\bfseries 50}, 23-47.
\bibitem{GuoZhang:2013} {\sc Guo, X. and Zhang, Y.} (2013). Generalized discounted continuous-time Markov decision
processes. arXiv:1304.3314.

\bibitem{Hernandez-Lerma:1996} {\sc Hern{\'a}ndez-Lerma, O. and Lasserre, J.} (1996). {\em Discrete-time {M}arkov Control Processes}. Springer, New York.
\bibitem{Hernandez-Lerma:2000} {\sc Hern{\'a}ndez-Lerma, O. and Lasserre, J.} (2000). {Fatou's lemma and Lebesgue's convergence theorem for measures}. {\it J. Appl. Math. Stoch. Anal.} {\bfseries 13}, 137-146.

\bibitem{JaskiewiczN:2006JMAA} {\sc Ja\'{s}kiewicz, A. and Nowak,
A.} (2006). On the optimality equation for average cost Markov
cntrol processes with Feller transition probabiliries. {\it J.
Math. Anal. Appl.} {\bfseries 316}, 495-509.

\bibitem{JaskiewiczN:2006OL} {\sc Ja\'{s}kiewicz, A. and Nowak,
A.} (2006). Optimality in Feller semi-Markov control processes.
{\it Oper. Res. Lett.} {\bfseries 34}, 713-718.
\bibitem{Jaskiewicz:2009} {\sc Ja\'{s}kiewicz, A.} (2009). Zero-sum ergodic semi-Markov games with weakly continuous transition probabilities.
{\it J. Optim. Theory Appl.} {\bfseries 141}, 321-347.

\bibitem{Kitaev:1986} {\sc Kitaev, M.} (1986). Semi-Markov and jump Markov controlled models: average cost criterion. {\it Theory. Probab. Appl.} {\bfseries 30}, 272-288.
\bibitem{Kitaev:1995} {\sc Kitaev, M. and Rykov, V.} (1995). {\em Controlled Queueing Systems}. CRC Press, Boca Raton.
\bibitem{Kuznetsov:1981} {\sc Kuznetsov, S.} (1981). Any Markov
process in a Borel space has a transition function. {\it Theory.
Probab. Appl.} {\bfseries 25}, 384-388.


\bibitem{ABPZY:20102} {\sc Piunovskiy, A. and Zhang, Y.} (2011). Discounted continuous-time Markov decision processes with unbounded rates: the convex analytic approach. {\it SIAM J. Control Optim.} {\bfseries 49}, 2032-2061.

\bibitem{PiunovskiyZY:2012} {\sc Piunovskiy, A. and Zhang, Y.} (2012). The transformation method for continuous-time Markov decision processes.  {\it J. Optim. Theory Appl.}, {\bfseries 154}, 691-712.

\bibitem{PrietoOHL:2012book} {\sc Prieto-Rumeau, T. and  Hern{\'a}ndez-Lerma,
O.} (2012). {\it Selected Topics in Continuous-time Controlled
Markov Chains and Markov Games}. Imperial College Press, London.
\bibitem{Puterman:1994} {\sc Puterman, M.} (1994). {\it {M}arkov {D}ecision {P}rocesses: {D}iscrete {S}tochastic {D}ynamic {P}rogramming.} Wiley, New York.


\bibitem{Zhu:2008}{\sc Zhu, Q.} (2008). Average optimality for continuous-time Markov decision processes with a policy iteration
approach. {\it J. Math. Anal. Appl.} {\bfseries 339}, 691-704.

\end{thebibliography}
\end{document}